\documentclass[12pt]{article}

\usepackage{theorem}
\usepackage{amssymb}
\usepackage{amsmath}
\usepackage{amsfonts}
\usepackage{times}
\usepackage{graphicx}

\def \RR {\mathbb R}

\def \EE {\mathbb E}

\def \vphi {\varphi}

\newtheorem{theorem}{Theorem}[section]
\newtheorem{inequality}{Inequality}[section]
\newtheorem{lemma}[theorem]{Lemma}
\newtheorem{question}[theorem]{Question}

\newtheorem{proposition}[theorem]{Proposition}
\newtheorem{corollary}[theorem]{Corollary}
 {\theorembodyfont{\rmfamily}}
\newtheorem{definition}[theorem]{Definition}

\def\myffrac#1#2 in #3{\raise 2.6pt\hbox{$#3 #1$}\mkern-1.5mu\raise 0.8pt\hbox{$
#3/$}\mkern-1.1mu\lower 1.5pt\hbox{$#3 #2$}}

\begin{document}

\title{Approximately gaussian marginals and the hyperplane conjecture}

\author{R. Eldan $\ $ and $\ $ B. Klartag\thanks{The authors were supported in part
by the Israel Science Foundation and by a Marie Curie Reintegration
Grant from the Commission of the European Communities.}}
\date{}
 \maketitle

\abstract{We discuss connections between certain well-known open
problems  related to the uniform measure on a high-dimensional
convex body. In particular, we show that the ``thin shell
conjecture'' implies the ``hyperplane conjecture''. This extends a
 result by K. Ball, according to which the stronger ``spectral gap
 conjecture'' implies the ``hyperplane conjecture''.}

\section{Introduction}
\label{sec1}

Little is currently known about the uniform measure on a general
high-dimensional convex body. Many aspects of the Euclidean ball or the unit cube
are easy to analyze, yet it is difficult to answer even some of
the simplest questions regarding arbitrary convex bodies, lacking
symmetries and structure. For example,

\begin{question}
Is there a universal constant $c > 0$ such that for any dimension $n$
 and a convex body $K \subset \RR^n$ with
$Vol_n(K) = 1$, there exists a hyperplane $H \subset \RR^n$ for
which $Vol_{n-1}(K \cap H) > c$? \label{slicing}
\end{question}

Here, of course, $Vol_k$ stands for $k$-dimensional volume. A convex
body is a bounded, open convex set. Question \ref{slicing} is
referred to as the ``slicing problem'' or the ``hyperplane
conjecture'', and was raised by Bourgain \cite{bou1, bou2} in
relation to the maximal function in high dimensions. It was
demonstrated  by Ball \cite{Ball} that Question \ref{slicing} and
similar questions are most naturally formulated in the broader class
of logarithmically concave densities.

\medskip A probability density $\rho: \RR^n \rightarrow [0, \infty)$
is called {\it log-concave} if it takes the form $\rho = \exp(-H)$
for a convex function $H: \RR^n \rightarrow \RR \cup \{ \infty \}$.
A probability measure is log-concave if it has a log-concave density.
The uniform probability measure on a convex body is an example of a
log-concave probability measure, as well as the standard gaussian
measure on $\RR^n$. A log-concave probability density decays
exponentially at infinity (e.g., \cite[Lemma 2.1]{K_psitwo}), and thus has moments of all orders. For a
probability measure $\mu$ on $\RR^n$ with finite second moments, we
consider its barycenter $b(\mu) \in \RR^n$ and covariance matrix $Cov(\mu)$
defined by
$$ b(\mu) = \int_{\RR^n} x d \mu(x), \ \ \ \ \ \ Cov(\mu) = \int_{\RR^n} (x - b(\mu)) \otimes (x - b(\mu)) d \mu(x) $$
where for $x \in \RR^n$ we write $x \otimes x$ for the $n\times n$
matrix $(x_i x_j)_{i,j=1,\ldots,n}$. A log-concave probability
measure $\mu$ on $\RR^n$ is {\it isotropic} if its barycenter lies
at the origin and its covariance matrix is the identity matrix. For
an isotropic, log-concave probability measure $\mu$ on $\RR^n$ we
denote
$$ L_{\mu} = L_f = f(0)^{1/n} $$
where $f$ is the log-concave density of $\mu$. It is well-known
(see, e.g., \cite[Lemma 3.1]{K_psitwo}) that $L_f > c$, for some
universal constant $c > 0$. Define
$$ L_n = \sup_{\mu} L_{\mu} $$
where the supremum runs over all isotropic, log-concave probability
measures $\mu$ on $\RR^n$. As follows from the works of Ball \cite{Ball}, Bourgain \cite{bou1}, Fradelizi \cite{fra}, Hensley \cite{Hensley} and Milman and Pajor \cite{MP}, Question \ref{slicing} is directly equivalent to the following:

\begin{question} Is it true that $\sup_n L_{n} < \infty$? \label{slicing2}
\end{question}

See also Milman and Pajor \cite{MP} and the second author's paper \cite{K_quarter} for a survey of results revolving around this question.  For a convex body $K \subset \RR^n$ we write
 $\mu_K$ for the uniform probability measure on $K$. A convex body $K \subset \RR^n$ is centrally-symmetric if $K = -K$.
 It is known that
\begin{equation}
 L_n \leq C \sup_{K \subset \RR^n} L_{\mu_K} \label{eq_1316}
\end{equation}
  where the supremum runs over all centrally-symmetric convex
bodies $K \subset \RR^n$ for which $\mu_K$ is isotropic, and $C > 0$
is a universal constant. Indeed, the reduction from log-concave
distributions to convex bodies was proven by Ball \cite{Ball} (see
\cite{K_quarter} for the straightforward generalization to the non-symmetric
case), and the reduction from general convex bodies to
centrally-symmetric ones was outlined, e.g., in the last paragraph
of \cite{K_funct}. The best estimate known to date is $L_{n} < C
n^{1/4}$ for a universal constant $C > 0$ (see \cite{K_quarter}),
which slightly sharpens an earlier estimate by Bourgain \cite{bou3,
bou_psi2, dar}.

\medskip Our goal in this note is to establish a connection between the
slicing problem and another  open problem in high-dimensional
convex geometry. Write $| \cdot |$ for the standard Euclidean norm
in $\RR^n$, and denote by $x \cdot y$ the scalar product of $x, y
\in \RR^n$.  We say that a random vector $X$ in $\RR^n$ is isotropic
and log-concave if it is distributed according to an isotropic,
log-concave probability measure.  Let $\sigma_n \geq 0$ satisfy
\begin{equation}
 \sigma_n^2 = \sup_X  \EE (|X| - \sqrt{n})^2 \label{eq_1356}
\end{equation}
  where the supremum
runs over all isotropic, log-concave random vectors $X$ in $\RR^n$.
The parameter $\sigma_n$ measures the width of the ``thin spherical
shell'' of radius $\sqrt{n}$ in which most of the mass of $X$ is
located. See (\ref{eq_958}) below for another definition of
$\sigma_n$, equivalent up to a universal constant, which is perhaps
more common in the literature. It is known that $\sigma_n \leq C
n^{0.41}$ where $C > 0$ is a universal constant (see
\cite{K_power}), and it is suggested in the works of Anttila, Ball
and Perissinaki \cite{ABP} and of Bobkov and Koldobsky \cite{BK}
that perhaps
\begin{equation} \sigma_n \leq C \label{eq_1344}
\end{equation} for a universal constant $C
> 0$. Again, up to a universal constant, one may restrict
attention in (\ref{eq_1356}) to random vectors that are distributed
uniformly in centrally-symmetric convex bodies. This essentially
follows from the same technique as in the case of the parameter
$L_n$ mentioned above.

\medskip The importance of the parameter $\sigma_n$ stems from the central
limit theorem for convex bodies \cite{K1}. This theorem
 asserts that most of the one-dimensional marginals of an isotropic,
 log-concave random vector are
approximately gaussian. The Kolmogorov distance to the standard
gaussian distribution of a typical marginal has roughly
the order of magnitude of $\sigma_n / \sqrt{n}$. Therefore, the conjectured
bound (\ref{eq_1344}) actually concerns the quality of
the gaussian approximation to the marginals of high-dimensional
log-concave measures. Our main result reads as follows:

\begin{inequality} For any $n \geq 1$,
\begin{equation}
 L_n \leq C \sigma_n
 \label{desired}
 \end{equation}
where $C > 0$ is a universal constant. \label{main_thm}
\end{inequality}

Inequality \ref{main_thm} states, in particular, that an affirmative
answer to the slicing problem follows from the {\it thin shell
conjecture} (\ref{eq_1344}). This sharpens a result announced by
Ball \cite{ball}, according to which a positive answer
 to the slicing problem is implied by the  stronger
  conjecture suggested by Kannan, Lov\'asz and Simonovits \cite{KLS}. The
  quick argument leading from the latter conjecture to (\ref{eq_1344}) is explained
 in Bobkov and Koldobsky \cite{BK}.  Write $S^{n-1} = \{ x \in \RR^n ;
 |x| = 1 \}$ for the unit sphere, and denote
   $$ \underline{\sigma}_n = \frac{1}{\sqrt{n}} \sup_{X} \left| \EE X |X|^2 \right| = \frac{1}{\sqrt{n}} \sup_X
   \sup_{\theta \in S^{n-1}} \EE (X \cdot \theta) |X|^2, $$
where the supremum runs over all isotropic, log-concave random vectors $X$ in $\RR^n$.

\begin{lemma} For any $n \geq 1$,
\begin{equation}
\sigma_n^2 \leq \frac{1}{n} \sup_X \EE (|X|^2 - n)^2 \leq C
\sigma_n^2, \label{eq_958}
\end{equation}
 where the supremum runs over all isotropic, log-concave random vectors $X$ in
$\RR^n$. Furthermore,
$$ 1 \leq \underline{\sigma}_n \leq C \sigma_n \leq \tilde{C} n^{0.41}. $$
Here, $C, \tilde{C} > 0$ are universal constants.
 \label{trivial}
\end{lemma}

Inequality \ref{main_thm} may be sharpened, in view of Lemma
\ref{trivial}, to the bound $$ L_n \leq C \underline{\sigma}_n, $$ for a
universal constant $C > 0$. This is explained in the proof of Inequality \ref{main_thm} in Section \ref{sec3}. Our argument involves a certain
Riemannian structure, which is presented in Section \ref{sec2}.

\medskip As the reader has probably already guessed, we use the letters $c,
\tilde{c}, c^{\prime}, C, \tilde{C}, C^{\prime}$ to denote positive
universal constants, whose value is not necessarily the same in
different appearances. Further notation and facts to be used
throughout the text: The support $Supp(\mu)$ of a Borel measure
$\mu$ on $\RR^n$ is the minimal closed set of full measure. When
$\mu$ is log-concave, its support is a convex set. For a Borel
measure $\mu$ on $\RR^n$ and a Borel map $T: \RR^n \rightarrow
\RR^k$ we define the push-forward of $\mu$ under $T$ to be the
measure $\nu = T_*(\mu)$ on $\RR^k$ with
$$ \nu(A) = \mu(T^{-1}(A)) \ \ \ \ \ \text{for any Borel set} \ A \subset \RR^k. $$
 Note that for any log-concave probability measure $\mu$ on
$\RR^n$, there exists an invertible affine map $T: \RR^n \rightarrow
\RR^n$ such that $T_*(\mu)$ is isotropic. When $T$ is a linear
function and $k < n$, we say that $T_*(\mu)$ is a marginal of $\mu$.
The Pr\'ekopa-Leindler inequality implies that any marginal of a
log-concave probability measure is itself a log-concave probability
measure. The Euclidean unit ball is denoted by $B_2^n = \{ x \in
\RR^n ; |x| \leq 1 \}$, and its volume satisfies
$$ \frac{c}{\sqrt{n}} \leq Vol_n(B_2^n)^{1/n} \leq \frac{C}{\sqrt{n}}.
$$
We write $\nabla \vphi$ for the gradient of the function $\vphi$,
and $\nabla^2 \vphi$ for the hessian matrix.
For $\theta \in S^{n-1}$ we write $\partial_{\theta}$ for
differentiation in direction $\theta$, and $\partial_{\theta \theta}(\vphi) = \partial_{\theta}(\partial_{\theta} \vphi)$.

\medskip {\emph{ Acknowledgements.}} We would like to thank Daniel Dadush, Vitali Milman,
Leonid Polterovich, Misha Sodin and Boris Tsirelson for interesting discussions related to this work,
and to Shahar Mendelson for pointing out that there is a difference
between extremal points and exposed points.

\section{A Riemannian metric associated with a convex body}
\label{sec2}

The main mathematical idea presented in this note is a certain
Riemannian metric associated with a convex body $K \subset \RR^n$.
Our construction is affinely invariant: We actually associate a
Riemannian metric with any affine equivalence class of convex bodies
(two convex bodies in $\RR^n$ are affinely equivalent if there
exists an invertible affine transformation  that maps one to the
other. Thus, all ellipsoids are affinely equivalent).

\medskip
Begin by recalling the technique from \cite{K_quarter}. Suppose
that $\mu$ is a compactly-supported Borel probability measure on $\RR^n$
whose support is not contained in a hyperplane. Denote by $K \subset
\RR^n$ the interior of the convex hull of $Supp(\mu)$, so $K$ is a convex body. The
{\it logarithmic Laplace transform} of $\mu$ is
\begin{equation}
 \Lambda(\xi) = \Lambda_{\mu}(\xi) = \log \int_{\RR^n} \exp(\xi
\cdot x) d \mu(x) \ \ \ \ \ \ \ \ \ \ \ (\xi \in \RR^n).
\label{eq_1146}
\end{equation}
 The
function $\Lambda$ is strictly convex and $C^{\infty}$-smooth on
$\RR^n$. For $\xi \in \RR^n$ let $\mu_{\xi}$ be the probability
measure on $\RR^n$ for which the density $d \mu_{\xi} / d \mu$ is
proportional to $x \mapsto \exp( \xi \cdot x )$. Differentiating
under the integral sign, we see that
$$ \nabla \Lambda(\xi) = b(\mu_\xi), \ \ \ \ \ \
\nabla^2 \Lambda(\xi) = Cov(\mu_\xi) \ \ \ \ \ \ \ (\xi \in \RR^n),
$$
where $b(\mu_\xi)$ is the barycenter of the probability measure
$\mu_{\xi}$ and $Cov(\mu_{\xi})$ is the covariance matrix. We
learned the following lemma from Gromov's work \cite{g}. A proof is
provided for the reader's convenience.

\begin{lemma} In the above notation,
$$ \int_{\RR^n} \det \nabla^2 \Lambda(\xi) d \xi =
Vol_n(K). $$ \label{lem_957}
\end{lemma}

\emph{Proof:} It is well-known that the open set $\nabla \Lambda(\RR^n) =
\{ \nabla \Lambda(\xi) ; \xi \in \RR^n \}$ is convex, and
that the map $\xi \mapsto \nabla \Lambda(\xi)$ is one-to-one
 (see, e.g., Rockafellar \cite[Theorem 26.5]{R}). Denote by $\overline{K}$
the closure of $K$. Then,
 \begin{equation}
 \nabla \Lambda(\RR^n) \subseteq \overline{K} \label{eq_1009}
 \end{equation}
 since for any $\xi \in \RR^n$, the point $\nabla \Lambda(\xi) \in \RR^n$ is
 the barycenter of a certain probability measure supported on the compact, convex set
 $\overline{K}$. Next we show that  $\overline{\nabla \Lambda(\RR^n)}$
 contains all of the exposed points of $Supp(\mu)$. Let $x_0 \in Supp(\mu)$ be an exposed
 point, i.e., there exists $\xi \in \RR^n$ such that
 \begin{equation}
  \xi \cdot x_0 > \xi \cdot x \ \ \ \ \ \ \text{for all} \ x_0 \neq x \in
 Supp(\mu). \label{eq_1028}
\end{equation}
We claim that \begin{equation} \lim_{r \rightarrow \infty} \nabla
\Lambda(r \xi) = x_0. \label{eq_1024}
\end{equation}
Indeed, (\ref{eq_1024}) follows from (\ref{eq_1028}) and from the
fact that $x_0$ belongs to the support of $\mu$: The measure $\mu_{r
\xi}$ converges weakly to the delta measure $\delta_{x_0}$ as $r
\rightarrow \infty$, hence the barycenter of $\mu_{r \xi}$ tends to
$x_0$. Therefore $x_0 \in \overline{\nabla \Lambda(\RR^n)}$. Any
exposed point of $\overline{K}$ is an exposed point of $Supp(\mu)$,
and we conclude that all of the exposed points of $\overline{K}$ are
contained in $\overline{\nabla \Lambda(\RR^n)}$. From Straszewicz's
theorem
 (see, e.g., Schneider \cite[Theorem 1.4.7]{S}) and from (\ref{eq_1009}) we deduce that
 $$ \overline{K} = \overline{\nabla \Lambda(\RR^n)}. $$
The set  $\nabla \Lambda(\RR^n)$ is open and convex, hence
necessarily $\nabla \Lambda(\RR^n) = K$. Since $\Lambda$ is
strictly-convex, its hessian is positive-definite everywhere, and
according to  the change of variables formula,
$$ Vol_n(K) = Vol_n \left( \nabla \Lambda(\RR^n) \right)
= \int_{\RR^n} \det \nabla^2 \Lambda(\xi) d \xi. $$ \hfill $\square$

\medskip Recall that $\mu$ is any compactly-supported probability measure on
$\RR^n$ whose support is not contained in a hyperplane.  For each
$\xi \in \RR^n$ the hessian matrix $\nabla^2 \Lambda(\xi) =
Cov(\mu_{\xi})$ is positive definite. For $\xi \in \RR^n$  set
\begin{equation}
 g(\xi)(u, v) = g_{\mu}(\xi)(u,v) = Cov(\mu_\xi)u \cdot v  \ \ \ \ \ \ \ \ \ \ (u,v \in \RR^n). \label{eq_2118}
\end{equation}
Then $g_{\mu}(\xi)$ is a
positive-definite bilinear form for any $\xi \in \RR^n$, and thus $g_{\mu}$ is a Riemannian metric
on $\RR^n$. We also set
\begin{equation}
 \Psi_{\mu}(\xi) = \log \frac{\det \nabla^2 \Lambda(\xi)}{\det
\nabla^2 \Lambda(0)} = \log \frac{\det Cov(\mu_\xi)}{\det Cov(\mu)}
\ \ \ \ \ \ \ \ \ \ \ \ \ \ \ (\xi \in \RR^n). \label{eq_1733}
\end{equation}
 We say
that $X_{\mu} = (\RR^n, g_{\mu}, \Psi_{\mu}, 0)$ is the ``Riemannian
package associated with the measure $\mu$''.

\begin{definition} A ``Riemannian package of dimension $n$'' is a quadruple $X =
(U, g, \Psi, x_0)$ where $U \subset \RR^n$ is an open set, $g$ is a
 Riemannian metric on $U$,  $x_0 \in
 U$ and $\Psi: U \rightarrow \RR$ is a function with $\Psi(x_0) =
 0$.
\end{definition}

Suppose $X = (U, g, \Psi, x_0)$ and $Y = (V, h, \Phi, y_0)$ are
Riemannian packages. A map $\vphi: U \rightarrow V$ is an
isomorphism of $X$ and $Y$ if the following conditions hold:
\begin{enumerate}
\item $\vphi$ is a Riemannian isometry between the
Riemannian manifolds $(U, g)$ and $(V, h)$.
\item $\vphi(x_0) = y_0$.
\item $\Phi( \vphi(x) ) = \Psi(x)$ for any $x \in U$.
\end{enumerate}
When such an isomorphism exists we say that $X$ and $Y$ are isomorphic, and we write $X
\cong Y$.

\medskip Let us describe an additional construction of the same Riemannian package
associated with $\mu$, a construction which is dual to the one
mentioned above.
 Consider the
Legendre transform
$$ \Lambda^*(x) = \sup_{\xi \in \RR^n} \left[ \xi \cdot x -
\Lambda(\xi) \right] \ \ \ \ \ \ \ \ \ (x \in K).
$$
Then $\Lambda^*: K \rightarrow \RR$ is a
strictly-convex $C^{\infty}$-function, and $\nabla \Lambda^*:
K \rightarrow \RR^n$ is the inverse map of $\nabla
\Lambda:\RR^n \rightarrow K$ (see Rockafellar \cite[Chapter
V]{R}). Define
$$ \Phi_{\mu}(x) = \log \frac{\det \nabla^2 \Lambda^*(b(\mu))}{\det \nabla^2 \Lambda^*(x)}
 \ \ \ \ \ \ \ \ \ \ (x \in K), $$ and for $x \in K$ set
$$ h(x)(u,v) = h_{\mu}(x)(u,v) = \left[ \nabla^2 \Lambda^*\right] (x)u \cdot v \ \ \ \ \ \ \ \ \ \ \ \ (u, v \in \RR^n). $$
Then $h_{\mu}$ is a Riemannian metric on $K$. Note the identity
$$ \left[ \nabla^2
\Lambda (\xi) \right]^{-1} = \left[ \nabla^2 \Lambda^* \right]( \nabla \Lambda(\xi) )
\ \ \ \ \ \ \ \ \ \  \ \ \ \ (\xi \in \RR^n).
$$
Using this identity, it is a simple exercise to verify that the
Riemannian package $\tilde{X}_{\mu} = (K, h_{\mu},
\Phi_{\mu}, b(\mu))$ is isomorphic to the Riemannian package
$X_{\mu} = (\RR^n, g_{\mu}, \Psi_{\mu}, 0)$ described earlier,
with $x = \nabla \Lambda(\xi)$ being the isomorphism.

\medskip The constructions $X_{\mu}$ and $\tilde{X}_{\mu}$ are
equivalent, and each has advantages over the other. It seems that
$X_{\mu}$ is preferable when carrying out computations, as the
notation is usually less heavy in this case. On the other hand, the
definition $\tilde{X}_{\mu}$ is perhaps easier to visualize: Suppose
that $\mu$ is the uniform probability measure on $K$. In this case
$\tilde{X}_{\mu}$ equips the convex body $K$ itself with a
Riemannian structure. One is thus tempted to imagine, for instance, how geodesics
look on $K$, and what is a Brownian motion in the body $K$ with respect to this metric. The following
lemma shows that this Riemannian structure on $K$ is invariant under
linear transformations.

\begin{lemma} Suppose $\mu$ and $\nu$ are compactly-supported probability measures on
$\RR^n$ whose support is not contained in a hyperplane. Assume that
there exists a linear map $T: \RR^n \rightarrow \RR^n$ such that
$$ \nu = T_*(\mu). $$
Then $X_{\mu} \cong X_{\nu}$. \label{trivial_lem}
\end{lemma}

\emph{Proof:} It is straightforward to check that the linear map
$T^t$ (the transposed matrix) is the required isometry between the
Riemannian manifolds $(\RR^n, g_{\nu})$ and $(\RR^n, g_{\mu})$.
However, perhaps a better way to understand this isomorphism, is to
note that the construction of $X_{\mu}$ may be carried out in a more
abstract fashion: Suppose that $V$ is an $n$-dimensional linear
space, denote by $V^*$ the dual space, and let $\mu$ be a
compactly-supported Borel probability measure on $V$ whose support
is not contained in a proper affine subspace of $V$. The logarithmic
Laplace transform $\Lambda: V^* \rightarrow \RR$ is well-defined, as
is the family of probability measures $\mu_{\xi} \ (\xi \in
V^*)$ on the space $V$. For a point $\xi \in V^*$ and two tangent
vectors $\eta, \zeta \in T_{\xi} V^* \equiv V^*$, set
\begin{equation}
 g_{\xi}(\eta, \zeta) = \int_{V} \eta(x) \zeta(x) d \mu_{\xi}(x) -
\left(  \int_{V} \eta(x)  d \mu_{\xi}(x) \right) \left(  \int_{V}
\zeta(x)  d \mu_{\xi}(x) \right). \label{eq_1602}
\end{equation}
 A moment of reflection reveals
that the definition (\ref{eq_1602}) of the positive-definite
bilinear form $g_{\xi}$ is equivalent to the definition
(\ref{eq_2118}) given above. Additionally, there exists a linear
operator $A_{\xi}: V^* \rightarrow V^*$, which is self-adjoint and
positive-definite with respect to the bilinear form $g_0$, that
satisfies
$$ g_{\xi}(\eta, \zeta) = g_0(A_{\xi} \eta, \zeta) \ \ \ \ \ \ \ \
\text{for all} \ \eta, \zeta \in V^*. $$ Hence we may define
$\Psi(\xi) = \log \det A_{\xi}$, which coincides with the definition
(\ref{eq_1733}) of $\Psi_{\mu}$ above. Therefore, $X_{\mu} = (V^*,
g, \Psi, 0)$ is the Riemannian package associated with $\mu$. Back
to the lemma, we see that $X_{\mu}$ is constructed from exactly the
same data as $X_{\nu}$, hence they must be isomorphic. \hfill $\square$

\begin{corollary} Suppose $\mu$ and $\nu$ are compactly-supported probability measures on
$\RR^n$ whose support is not contained in a hyperplane. Assume that
there exists an affine map $T: \RR^n \rightarrow \RR^n$ such that
$$ \nu = T_*(\mu). $$
Then $X_{\mu} \cong X_{\nu}$. \label{trivial_cor}
\end{corollary}

\emph{Proof:} The only difference from Lemma \ref{trivial_lem} is
that the map $T$ is assumed to be affine, and not linear. It is clearly enough to
deal with the case where $T$ is a translation, i.e.,
$$ T(x) = x + x_0 \ \ \ \ \ \ \ \ \ \ \ \ \ (x \in \RR^n) $$
for a certain vector $x_0 \in \RR^n$. From the definition
(\ref{eq_1146}) we see that
$$ \Lambda_{\nu}(\xi) = \xi \cdot x_0 + \Lambda_{\mu}(\xi) \ \ \ \ \
\ \ \ \ (\xi \in \RR^n). $$ Adding a linear functional does
not influence second derivatives, hence $g_{\mu} = g_{\nu}$ and also
$\Psi_{\mu} = \Psi_{\nu}$. Therefore $X_{\mu} = (\RR^n, g_{\mu},
\Psi_{\mu}, 0)$ is trivially isomorphic to $X_{\nu} = (\RR^n,
g_{\nu}, \Psi_{\nu}, 0)$. \hfill $\square$

\medskip An $n$-dimensional Riemannian package
 is of ``log-concave type'' if it is isomorphic to the
Riemannian package $X_{\mu}$ associated with a compactly-supported,
log-concave probability measure $\mu$ on $\RR^n$. Note that according to our  terminology, a log-concave probability measure is absolutely-continuous with respect to the Lebesgue measure on $\RR^n$, hence its support is never contained in a hyperplane.

\begin{lemma} Suppose $X = (U, g, \Psi, \xi_0)$ is an $n$-dimensional
Riemannian package of log-concave type. Let $\xi_1 \in U$.
Denote \begin{equation}  \tilde{\Psi}(\xi) = \Psi(\xi) - \Psi(\xi_1)
\ \ \ \ \ \ \ \ \ \ \ (\xi \in U). \label{eq_1136} \end{equation}
Then also $Y = (U, g, \tilde{\Psi}, \xi_1)$ is an $n$-dimensional
Riemannian package of log-concave type. \label{lem_1426}
\end{lemma}

\emph{Proof:} Let $\mu$ be a compactly-supported log-concave
probability measure on $\RR^n$ whose associated Riemannian package
$X_{\mu} = (\RR^n, g_\mu, \Psi_\mu, 0)$ is isomorphic to $X$. Thanks
to the isomorphism, we may identify $\xi_1$ with a certain point in
$\RR^n$, which will still be denoted by $\xi_1$ (with a slight abuse
of notation). We now interpret the definition (\ref{eq_1136}) as
$$ \tilde{\Psi}(\xi) = \Psi(\xi) -
\Psi(\xi_1) \ \ \ \ \ \ \ \ \ \ \ (\xi \in \RR^n). $$ In order to
prove the lemma, we need to demonstrate that
\begin{equation}
 Y = (\RR^n, g_{\mu}, \tilde{\Psi}, \xi_1) \label{eq_1138}
 \end{equation} is of log-concave
type.
 Recall that $\mu_{\xi_1}$ is the
compactly-supported probability measure on $\RR^n$ whose density
with respect to $\mu$ is proportional to $x \mapsto \exp( \xi_1 \cdot
x )$. A crucial observation is that $\mu_{\xi_1}$ is log-concave. Set $\nu = \mu_{\xi_1}$,
and note the relation
\begin{equation}
\Lambda_{\nu}(\xi) = \Lambda_{\mu}(\xi + \xi_1) -
\Lambda_{\mu}(\xi_1) \ \ \ \ \ \ \ \ \ \ \ \ \ \ \ \ (\xi \in \RR^n).
\label{eq_1148}
\end{equation}
It
suffices to show that the Riemannian package $Y$ in (\ref{eq_1138})
is isomorphic to $X_{\nu} = (\RR^n, g_\nu, \Psi_\nu, 0)$. We claim that an
isomorphism $\vphi$ between $X_{\nu}$ and $Y$ is simply the
translation $$ \vphi(\xi) = \xi + \xi_1 \ \ \ \ \ \ \ \ \ \ \
\ (\xi \in \RR^n). $$ In order to see that $\vphi$ is indeed an
isomorphism, note that (\ref{eq_1148}) yields
\begin{equation} \nabla^2 \Lambda_\nu(\xi) = \nabla^2
\Lambda_{\mu}(\xi + \xi_1) \ \ \ \ \ \ \ \ \ \ \ \ \ \ \ \ (\xi \in
\RR^n), \label{eq_1157}
\end{equation} hence $\vphi$ is a Riemannian isometry between
$(\RR^n, g_{\nu})$ and $(\RR^n, g_{\mu})$, with $\vphi(0)
= \xi_1$. The relation (\ref{eq_1157}) implies that
$\tilde{\Psi}(\vphi(\xi)) = \Psi_{\nu}(\xi)$ for all $\xi \in \RR^n$. Hence
$\vphi$ is an isomorphism between Riemannain packages, and the lemma
is proven. \hfill $\square$

\medskip {\emph {Remark.}} When $\mu$ is a product measure on $\RR^n$ (such as
the uniform probability measure on the cube, or the gaussian
measure), straightforward computations of curvature show that the
manifold $(\RR^n, g_{\mu})$ is flat (i.e., all sectional curvatures
vanish). We were not able to extract meaningful information from the
local structure of the Riemannian manifold $(\RR^n, g_{\mu})$ in the
general case.

\section{Inequalities}
\label{sec3}

\emph{Proof of Lemma \ref{trivial}:} First, note that for any random
vector $X$ in $\RR^n$ with finite fourth moments,
$$
 \EE (|X| - \sqrt{n})^2 \leq \frac{1}{n} \EE (|X| - \sqrt{n})^2
(|X| + \sqrt{n})^2 = \frac{1}{n} \EE (|X|^2 - n)^2. $$ This proves
the inequality on the left in (\ref{eq_958}).  Regarding the
 inequality on the right, we use the bound
\begin{equation}
 \EE |X|^4 1_{|X| > C \sqrt{n}} \leq C \exp
\left(-\sqrt{n} \right) \label{eq_1049} \end{equation}
which follows from Paouris theorem \cite{Pa}. Here $1_{|X| > C \sqrt{n}}$ is the random variable that equals one
when $|X| > C \sqrt{n}$ and vanishes otherwise. Apply
 again the identity $|X|^2 - n =
(|X| - \sqrt{n}) (|X| + \sqrt{n})$ to conclude that
\begin{eqnarray}
\lefteqn{  \nonumber \EE (|X|^2 - n)^2 = \EE (|X|^2 - n)^2 1_{|X|
\leq C \sqrt{n}} + \EE (|X|^2 - n)^2 1_{|X| > C \sqrt{n}} }
\\ & \leq & (C + 1)^2 n \EE
 (|X| - \sqrt{n})^2  + \EE |X|^4 1_{|X| > C \sqrt{n}},
 \label{eq_1048}
 \phantom{aaaaaaaaaaaaaa}
\end{eqnarray}
where $C \geq 1$ is the universal constant from (\ref{eq_1049}). A simple
computation shows that  $\sigma_n \geq \sqrt{2}$, as is
witnessed by the standard gaussian random vector in $\RR^n$, or by
the example in the next paragraph. Thus the inequality on the right
in (\ref{eq_958}) follows from (\ref{eq_1049}) and (\ref{eq_1048}).
Our proof of (\ref{eq_958}) utilized the deep Paouris theorem. Another
possibility could be to use \cite[Theorem 4.4]{K_power} or the
deviation inequalities for polynomials  proved first by Bourgain
\cite{bou3}.

\medskip In order to prove the second assertion in the lemma, observe that since $\EE X = 0$,
$$ \EE (X \cdot \theta) |X|^2 =
\EE (X \cdot \theta) (|X|^2 - n)  \leq  \sqrt{ \EE (X \cdot
\theta)^2 \EE (|X|^2 - n)^2 } \leq C \sqrt{n} \sigma_n,
$$
where we used the Cauchy-Schwartz inequality, the fact that $\EE (X
\cdot \theta)^2 = 1$ and (\ref{eq_958}). It remains to prove that
$\underline{\sigma}_n \geq 2$. To this end, consider the case where
$Y_1,\ldots,Y_n$ are independent random variables, all distributed
according to the density $t \mapsto e^{-I(t + 1)}$ on the real line,
where $I(a) = a$ for $a \geq 0$ and $I(a) = +\infty$ for $a < 0$.
Then $Y = (Y_1,\ldots,Y_n)$ is a random vector distributed according
to an isotropic, log-concave probability measure on $\RR^n$, and
$$ \EE \frac{\sum_{j=1}^n Y_j}{\sqrt{n}} |Y|^2 = 2 \sqrt{n}. $$
This completes the proof. \hfill $\square$

\medskip When $\vphi$ is a smooth real-valued function on a Riemannian manifold
$(M, g)$, we denote its gradient at the point $x_0 \in M$ by
$\nabla_g \vphi(x_0) \in T_{x_0}(M)$. Here $T_{x_0}(M)$ stands for the
tangent space to $M$ at the point $x_0$. The subscript $g$ in
$\nabla_g \vphi(x_0)$ means that the gradient is computed with
respect to the Riemannian metric $g$. The usual gradient of a
function $\vphi:\RR^n \rightarrow \RR$ at a point $x_0 \in \RR^n$
is denoted by $\nabla \vphi(x_0) \in \RR^n$, without any subscript.
The length of a tangent vector $v \in T_{x_0}(M)$ with respect to the metric $g$ is $|v|_g = \sqrt{g_{x_0}(v,v)}$.

\begin{lemma} Suppose $X = (U, g, \Psi, \xi_0)$ is an $n$-dimensional
Riemannian package of log-concave type. Then, for any $\xi \in U$,
$$ |\nabla_g \Psi(\xi)|_g \leq \sqrt{n} \underline{\sigma}_n. $$ \label{lem_1427}
\end{lemma}

\emph{Proof:} Suppose first that $\xi = \xi_0$. We  need to
establish the bound
\begin{equation}
 |\nabla_g \Psi(\xi_0)|_g \leq \sqrt{n} \underline{\sigma}_n \label{eq_1749}
 \end{equation}
for any log-concave package $X = (U, g, \Psi, \xi_0)$ of dimension
$n$. Any such package $X$ is isomorphic to $X_{\mu} = (\RR^n,
g_{\mu}, \Psi_{\mu}, 0)$ for a certain compactly-supported
log-concave probability measure $\mu$ on $\RR^n$. Furthermore,
according to Corollary \ref{trivial_cor}, we may apply an
appropriate affine map and assume that $\mu$ is isotropic. Thus our
goal is to prove that
\begin{equation}
 |\nabla_{g_{\mu}} \Psi_{\mu}(0)|_{g_\mu} \leq \sqrt{n} \underline{\sigma}_n. \label{eq_1451} \end{equation}
Since $\mu$ is
isotropic,  $\nabla^2 \Lambda_{\mu}(0) = Cov(\mu) = Id$, where $Id$ is
the identity matrix. Consequently, the desired bound (\ref{eq_1451}) is
equivalent to
$$
 |\nabla \Psi_{\mu}(0)| \leq \sqrt{n} \underline{\sigma}_n.
$$  Equivalently, we need to show that
$$
\left. \partial_{\theta} \log \frac{\det \nabla^2
\Lambda_{\mu}(\xi)} {\det \nabla^2 \Lambda_{\mu}(0)} \right|_{\xi =
0} \leq \sqrt{n} \underline{\sigma}_n \ \ \ \ \ \ \ \ \ \ \ \ \
\text{for all} \ \theta \in S^{n-1}. $$ A straightforward
computation shows that $\partial_{\theta} \log \det \nabla^2
\Lambda_{\mu}(\xi)$ equals the trace of the matrix $\left( \nabla^2
\Lambda_{\mu}(\xi) \right)^{-1} \nabla^2
\partial_{\theta} \Lambda_{\mu}(\xi)$. Since $\mu$ is isotropic,
$$ \left. \partial_{\theta} \log \frac{\det \nabla^2 \Lambda_{\mu}(\xi)}
{\det \nabla^2 \Lambda_{\mu}(0)} \right|_{\xi = 0} = \triangle
\partial_{\theta} \Lambda_{\mu}(0) = \int_{\RR^n} (x \cdot \theta)
|x|^2 d \mu(x) \leq \sqrt{n} \underline{\sigma}_n, $$ according to
the definition of $\underline{\sigma}_n$, where $\triangle$ stands
for the usual Laplacian in $\RR^n$. This completes the proof of
(\ref{eq_1749}). The lemma in thus proven in the special case where
$\xi = \xi_0$.

\medskip The general case follows from  Lemma \ref{lem_1426}:
When $\xi \neq \xi_0$, we may consider the log-concave Riemannian package $Y =
(U, g, \tilde{\Psi}, \xi)$, where $\tilde{\Psi}$ differs from $\Psi$
by an additive constant. Applying (\ref{eq_1749}) with the
log-concave package $Y$, we see that
$$ |\nabla_g \Psi(\xi)|_g = |\nabla_g \tilde{\Psi}(\xi)|_g \leq
\sqrt{n} \underline{\sigma}_n. $$ \hfill $\square$

\medskip The next lemma is a crude upper bound for the
Riemannian distance, valid for any Hessian metric (that is, a
Riemannian metric on $U \subset \RR^n$ induced by the
hessian of a convex function).

\begin{lemma} Let $\mu$ be a compactly-supported
probability measure on $\RR^n$ whose support is not contained in a
hyperplane. Denote by $\Lambda$ its logarithmic Laplace transform,
and let $X_{\mu} = (\RR^n, g_{\mu}, \Psi_{\mu}, 0)$ be the
associated Riemannian package. Then for any $\xi, \eta \in \RR^n$,
\begin{equation}
 d(\xi, \eta) \leq \sqrt{ \Lambda( 2 \xi - \eta) - \Lambda(\eta) - 2 \nabla \Lambda(\eta) \cdot (\xi - \eta)
 },
\label{eq_1505}
\end{equation}
where $d(\xi, \eta)$ is the Riemannian distance between $\xi$ and $\eta$, with respect to the Riemannian metric $g_{\mu}$.
In particular, when the barycenter of $\mu$ lies at the origin,
\begin{equation}
 d(\xi, 0) \leq \sqrt{\Lambda(2 \xi)}. \label{eq_1527}
\end{equation}
\label{lem_1150}
\end{lemma}

\emph{Proof:} The bound (\ref{eq_1505}) is obvious when $\xi = \eta$.
When $\xi \neq \eta$, we need to exhibit a path from $\eta$ to
$\xi$ whose Riemannian length is at most the expression on the right in (\ref{eq_1505}). Set
$\theta = (\xi - \eta) / |\xi - \eta|$ and $R = |\xi - \eta|$. Consider the interval
$$ \gamma(t) = \eta + t \theta \ \ \ \ \ \ \ \ \ \ \ \ \ (0 \leq t \leq R). $$
This path connects $\eta$ and $\xi$, and its Riemannian
length is
\begin{eqnarray*}
 \lefteqn{ \int_0^{R} \sqrt{ g_{\mu}(\gamma(t)) \left(
\theta, \theta \right) } dt = \int_0^{R} \sqrt{ [\partial_{\theta \theta} \Lambda]  (\eta + t \theta) }
dt} \\ & = & \int_0^{R} \sqrt{ \frac{d^2 \Lambda  (\eta + t \theta)}{dt^2}  }
dt \leq \sqrt { \int_0^{2R} (2 R - t) \frac{d^2 \Lambda  (\eta + t \theta)}{dt^2}  dt
\int_0^R \frac{dt}{2 R - t} },
\end{eqnarray*}
according to the Cauchy-Schwartz inequality. Clearly, $\int_0^R dt /
(2 R - t) = \log 2 \leq 1$. Regarding the other integral, recall
Taylor's formula with integral remainder:
$$ \int_0^{2R} (2 R - t) \frac{d^2 \Lambda  (\eta + t \theta)}{dt^2}  dt =
\Lambda(\eta + 2 R \theta) - \left[ \Lambda(\eta) + 2 R \theta \cdot
\nabla \Lambda(\eta) \right]. $$ The inequality (\ref{eq_1505}) is thus proven. Furthermore, $\Lambda(0)
= 0$, and when the barycenter of $\mu$ lies at the origin, also
$\nabla \Lambda(0) = 0$. Thus (\ref{eq_1527}) follows from
(\ref{eq_1505}). \hfill $\square$

\medskip The volume-radius of a convex body $K \subset \RR^n$ is $$ v.rad.(K) =
(Vol_n(K) / Vol_n(B_2^n))^{1/n}. $$ This is the radius of the
Euclidean ball that has exactly the same volume as $K$. When $E
\subseteq \RR^n$ is an affine subspace of dimension $\ell$ and $K
\subset E$ is a convex body, we interpret $v.rad.(K)$ as $(Vol(K) /
Vol(B_2^\ell))^{1/\ell}$. For a subspace $E \subset \RR^n$, denote
by $Proj_E: \RR^n \rightarrow E$ the orthogonal projection operator
onto $E$ in $\RR^n$.
 A Borel measure $\mu$ on $\RR^n$ is {\it even} or {\it centrally-symmetric} if $\mu(A) =
\mu(-A)$ for any measurable $A \subset \RR^n$.

\begin{lemma} Let $\mu$ be an even, isotropic, log-concave probability measure on $\RR^n$.
Let $1 \leq t \leq \sqrt{n}$ and denote by $B_t \subset \RR^n$ the
collection of all $\xi \in \RR^n$ with $d(0, \xi) \leq t$, where
$d(0, \xi)$ is as in Lemma \ref{lem_1150}. Then,
\begin{equation}
 Vol_n(B_t)^{1/n} \geq c \frac{t}{\sqrt{n}}, \label{eq_1511}
\end{equation}
where $c > 0$ is a universal constant. Here, as elsewhere, $Vol_n$
stands for the Lebesgue measure on $\RR^n$ (and not the Riemannian
volume). \label{cor_1055}
\end{lemma}

\emph{Proof:} It suffices to prove the lemma under the additional
assumption that $t$ is an integer. According to Lemma
\ref{lem_1150},
$$ K_t := \{ \xi \in \RR^n; \Lambda(2 \xi) \leq t^2 \} \subseteq B_t. $$
Let $E \subset \RR^n$ be any $t^2$-dimensional subspace, and denote
by $f_E: \RR^n \rightarrow [0, \infty)$ the density of the isotropic
probability measure $(Proj_E)_* \mu$. Then $f_E$ is a log-concave
function, according to the Pr\'ekopa-Leindler inequality, and $f_E$
is also an even function. According to the definition above,
$$
 f_E(0)^{1/t^2} = L_{f_E} \geq c. $$
Note that the restriction of $\Lambda$ to the subspace $E$ is the
logarithmic Laplace transform of $(Proj_E)_* \mu$. It is proven in
\cite[Lemma 2.8]{K_psitwo} that
\begin{equation}
 v.rad.(K_t \cap E) \geq c t f_E(0)^{1/t^2} \geq c^{\prime} t.
 \label{eq_1023}
 \end{equation}
 The bound (\ref{eq_1023}) holds for any subspace $E \subset \RR^n$
 of dimension $t^2$. From \cite[Corollary 3.1]{low_M} we deduce
 that
 $$ v.rad.(K_t) \geq \tilde{c} t. $$
Since $K_t \subseteq B_t$, the bound (\ref{eq_1511}) follows. \hfill
$\square$

\begin{lemma} Let $\mu$ be a compactly-supported, even, isotropic, log-concave
probability measure on $\RR^n$. Denote by $K$ the interior of the
support of $\mu$, a convex body in $\RR^n$. Then,
$$ Vol_n( K )^{1/n} \geq c  / \underline{\sigma}_n, $$
where $c > 0$ is a universal constant. \label{lem_1100}
\end{lemma}

\emph{Proof:} Set $t = \max \{ \sqrt{n} / \underline{\sigma}_n, 1 \}$.
Then $1 \leq t \leq \sqrt{n}$ and $\underline{\sigma}_n \leq C \sqrt{n}$, according to Lemma \ref{trivial}.
Recall the
definition of the set $B_t \subset \RR^n$ from Lemma
\ref{cor_1055}. Consider the Riemannian package $X_{\mu} = (\RR^n,
g_{\mu}, \Psi_{\mu}, 0)$ that is associated with the measure $\mu$.
According to Lemma \ref{lem_1427}, for any $\xi \in B_t$,
$$ \Psi_{\mu}(0) - \Psi_{\mu}(\xi) \leq \sqrt{n} \underline{\sigma}_n d (0, \xi) \leq t
\sqrt{n} \underline{\sigma}_n \leq Cn. $$ Since  $\Psi_{\mu}(\xi) = \log \det \nabla^2
\Lambda_{\mu}(\xi)$ and $\Psi_{\mu}(0) =
0$, then
$$ \det \nabla^2 \Lambda_{\mu}(\xi) \geq e^{-C n} \ \ \ \ \ \ \ \
\text{for any } \ \xi \in B_t. $$ From Lemma \ref{lem_957},
$$ Vol_n(K) = \int_{\RR^n} \det \nabla^2 \Lambda_{\mu}(\xi) d \xi
\geq \int_{B_t} \det \nabla^2 \Lambda_{\mu}(\xi) d \xi \geq e^{-C n}
Vol_n(B_t)
$$
as $\Lambda_{\mu}$ is convex and hence $\det \nabla^2
\Lambda_{\mu}(\xi) \geq 0$ for all $\xi$. Lemma \ref{cor_1055}
yields that
$$ Vol_n(K)^{1/n} \geq e^{-C} \left ( c \frac{t}{\sqrt{n}} \right ) \geq \frac{c^{\prime}}{\underline{\sigma}_n}. $$
The lemma is proven. \hfill $\square$

\medskip \emph{Proof of Inequality \ref{main_thm}}:
Let $K \subset \RR^n$ be a centrally-symmetric convex body
such that the uniform probability measure $\mu_K$ is isotropic. Then,
$$ L_{\mu_K} = \frac{1}{Vol_n(K)^{1/n}} \leq C \underline{\sigma}_n
$$
thanks to Lemma \ref{lem_1100}. In view of (\ref{eq_1316}), the
bound $L_n \leq C \underline{\sigma}_n$ is proven. The desired
inequality (\ref{desired}) now follows from Lemma \ref{trivial}. \hfill $\square$

\medskip The following proposition is not applied in this article. It is nevertheless
included as it may help understand the nature of the elusive
quantity $\left| \EE X |X|^2 \right|$ for an isotropic, log-concave
random vector $X$ in $\RR^n$.

\begin{proposition} Suppose $X$ is an isotropic random vector in $\RR^n$
with finite third moments. Then,
$$ \left| \EE X |X|^2 \right|^2 \leq C n^3 \int_{S^{n-1}} \left( \EE (X \cdot
\theta)^3 \right)^2 d \sigma_{n-1}(\theta) $$ where $\sigma_{n-1}$
is the uniform Lebesgue probability measure on the sphere $S^{n-1}$,
and $C > 0$ is a universal constant. \label{prop_1242}
\end{proposition}

\emph{Proof:} Denote $F(\theta) = \EE (X \cdot \theta)^3$ for
$\theta \in \RR^n$. Then $F(\theta)$ is a homogenous polynomial of
degree three, and its Laplacian is
$$ \triangle F(\theta) = 6 \EE (X \cdot \theta) |X|^2. $$
Denote $v = \EE X |X|^2 \in \RR^n$. The function
$$ \theta \mapsto F(\theta) - \frac{6}{2n+4} |\theta|^2 (\theta \cdot v) \ \ \ \ \ \ \ \ \ (\theta \in \RR^n) $$
is a homogenous, harmonic polynomial of degree three. In
other words, the restriction $F|_{S^{n-1}}$ decomposes
into spherical harmonics as
$$ F(\theta) = \frac{6}{2n+4} (\theta \cdot v) + \left( F(\theta) - \frac{6}{2n+1} (\theta \cdot
v) \right) \ \ \ \ \ \ \ \ \ \ (\theta \in S^{n-1}). $$ Since
spherical harmonics of different degrees are orthogonal to each
other,
$$ \int_{S^{n-1}} F^2 (\theta) d \sigma_{n-1}(\theta) \geq
\frac{36}{(2n+4)^2} \int_{S^{n-1}} (\theta \cdot v)^2  d
 \sigma_{n-1}(\theta) = \frac{36}{n (2n+4)^2} |v|^2. $$
\hfill $\square$

\medskip \emph{Remark.}
 According to Proposition \ref{prop_1242},
if we could show that $\left| \EE (X \cdot \theta)^3 \right| \leq C
/ n$ for a typical unit vector $\theta \in S^{n-1}$, we would obtain
a positive answer to Question \ref{slicing}. It is interesting to
note that the function
$$ F(\theta) = \EE |X \cdot \theta| \ \ \ \ \ \ \ \ \ \ \ \ \ \ (\theta \in S^{n-1}) $$
admits tight concentration bounds. For instance,
$$ \int_{S^{n-1}} (F(\theta) / E - 1)^2 d \sigma_{n-1}(\theta) \leq C /
n^2 $$ where $E = \int_{S^{n-1}} F(\theta) d \sigma_{n-1}(\theta)$,
whenever $X$ is distributed according to a suitably normalized
log-concave probability measure on $\RR^n$. The normalization we currently prefer here
is slightly different from the isotropic normalization. The
details will be explained elsewhere, as well as some relations to
the problem of stability in the Brunn-Minkowski inequality.

 {\small

\bigskip

{\small \noindent  School of Mathematical Sciences, Tel-Aviv
University, Tel-Aviv 69978, Israel

{\small \noindent \it e-mail address:} {\small
\verb"[roneneld,klartagb]@tau.ac.il"}


\end{document}